\documentclass[12pt]{article}
\usepackage{tikz}
\usetikzlibrary{arrows, decorations.markings, decorations.pathmorphing, backgrounds, positioning, fit, petri}
\usepackage{amsfonts}                 % Verwendung von \mathbb N
\usepackage{amssymb}
\usepackage[latin1]{inputenc}         % Tastatur
\usepackage{graphicx}
\newcommand{\cP}{{\cal P}}
\newcommand{\cQ}{{\cal Q}}
\newcommand{\cG}{{\cal G}}
\newcommand{\cT}{{\cal T}}
\newcommand{\cH}{{\cal H}}

\newcommand{\bewende}{\hspace*{\fill}$\Box $\\[0.5cm]} %Beweisende
\newcommand{\bewanf}{\noindent \textbf{Proof:}\quad }    %Beweisanfang
\newcommand{\map}[2]{$f\colon #1\to #2$}
\newcommand{\satzanf}{\begin{samepage}\begin{satz}}
\newcommand{\satzende}{\end{satz}\end{samepage}}

\newcommand{\propanf}{\begin{samepage}\begin{prop}}
\newcommand{\propende}{\end{prop}\end{samepage}}

\newcommand{\koranf}{\begin{samepage}\begin{kor}}
\newcommand{\korende}{\end{kor}\end{samepage}}
\newcommand{\lemanf}{\begin{samepage}\begin{lem}}
\newcommand{\lemende}{\end{lem}\end{samepage}}
\newcommand{\bspanf}{\begin{beispiel}}
\newcommand{\bspende}{\end{beispiel}}
\newcommand{\defanf}{\begin{definition}}
\newcommand{\defende}{\end{definition}}
\newtheorem{satz}{Theorem}[section]             %  anpassen f�r reportarticle
\newtheorem{lem}[satz]{Lemma}
\newtheorem{beispiel}[satz]{Example}
\newtheorem{definition}[satz]{Definition}
\newtheorem{kor}[satz]{Corollary}
\newtheorem{prop}[satz]{Proposition}

\newcommand{\be}{\begin{equation}}
\newcommand{\ee}{\end{equation}}
\newcommand{\bl}[1]{\begin{equation}\label{#1}}
\begin{document}
\title{Injective Labeled Oriented Trees are Aspherical}
\author{Jens Harlander and Stephan Rosebrock}
%\date{}                  %nur nehmen, wenn kein Datum darauf soll
\maketitle
\thispagestyle{empty}               % erste Seite ohne Kopf- und Fusszeile
\begin{abstract}
A labeled oriented tree is called {\it injective}, if each vertex occurs at most
once as an edge label. We show that injective labeled oriented trees are
aspherical. The proof uses a new relative 
asphericity test based on a lemma of Stallings.\\


\end{abstract}
%
%\tableofcontents                    % Inhaltsverzeichnis
%
\section{Introduction}

This article is concerned with the Whitehead conjecture, which states that a subcomplex of an aspherical 2-complex is aspherical. See Bogley \cite{Bogley} and Rosebrock \cite{Rosebrock2} for surveys. The conjecture originally arose in the context of knot theory.  The Wirtinger presentation of a knot gives rise to a 2-complex that is a subcomplex of a contractible 2-complex. Thus, an affirmative answer to the conjecture implies the asphericity of knot complements in the 3-sphere. Labeled oriented trees give rise to presentations that generalize Wirtinger presentations for knots, and presentations obtained from injective labeled oriented trees generalize Wirtinger presentations of alternating knots. Labeled oriented trees play a central role in understanding the Whitehead conjecture. Howie \cite{How83} showed that the finite case of the Whitehead conjecture reduces, up to the Andrews-Curtis conjecture, to the statement that presentations arising from labeled oriented trees are aspherical.\\

A {\em labeled oriented graph} (LOG) is an oriented graph $\cal G$ on vertices $\textbf{x}$ and edges ${\bf e}$, where
each oriented edge is labeled by a vertex. Associated with it is the {\em LOG-presentation} $P({\cal G})=\langle \textbf{x}\ |\ \{r_e\}_{e\in {\bf e}}\rangle$. If $e$ is an edge that starts at $x$, ends at $y$ and is labeled by $z$, then $r_e=xz(zy)^{-1}$. We also use the notation $e=[x,z,y]$.
A {\em LOG-complex} $K({\cal G})$ is the standard 2-complex associated with the LOG-presentation $P({\cal G})$, and a \emph{LOG-group} $G({\cal G})$, is the group defined by the LOG-presentation. We say a labeled oriented graph is {\em aspherical} if its associated LOG-complex is aspherical. A labeled oriented graph is called {\em injective} if each vertex occurs at most once as an edge label.
A {\em labeled oriented tree} (LOT) is a labeled oriented graph where the underlying graph is a tree. \\

The following is the main result of this article.

\satzanf\label{hauptsatz} Injective labeled oriented trees are aspherical.
\satzende

The theorem does not extend to labeled oriented graphs. The Wirtinger presentation $P$ read off an alternating knot diagram with $n$ crossings is the LOG-presentation of an injective labeled oriented circle $\mathcal C$ with $n$ edges. This labeled oriented circle $\mathcal C$ is not aspherical, because any one relator in $P$ is a consequence of the other relators.\\

We will need the following additional terminology on labeled oriented graphs. A {\em sub-LOG} of a labeled oriented graph $\cal G$ is a connected
subgraph $\cal H$ (containing at least one edge)  such that each edge label of $\cal H$ is a vertex of $\cal H$. A sub-LOG $\cal H$ is {\em proper} if it is not all of $\cal G$. 
A labeled oriented graph is called {\em compressed} if no edge is labeled with one of its vertices. It is called
{\em boundary reducible} if there is a boundary vertex that does not occur as edge label, and {\em boundary reduced} otherwise. 
A labeled oriented graph is called {\em interior reducible} if there is a vertex with two adjacent
edges with the same label that either point away or 
towards that vertex, and {\em interior reduced} otherwise. A labeled oriented graph which is boundary reduced, interior reduced and compressed is called {\em reduced}. 

Howie \cite{How85} observed that a labeled oriented tree $\cal G$ can be transformed into a reduced labeled oriented tree ${\cal G}_{red}$ so that $K(\cG)$ and $K(\cG_{red})$ have the same homotopy type. Here are some details on this transformation. If $\cG$ is not compressed then it contains an edge of the form $e=[a,a,b]$ or $e=[b,a,a]$. We remove the interior of $e$ and identify the vertices $a$ and $b$ to become one vertex $a$. Edges labeled with $b$ we relabel with $a$. This transforms $\cG$ to a labeled oriented tree $\cG'$ with fewer vertices.  
%On the presentation level $P(\cG)$ is changed into $P(\cG')$ by $Q^{**}$-transformations. 
On the 2-complex level the change amounts to a 3-deformation $h\colon K(\cG) \to K(\cG')$. If $\cG$ is not boundary reduced, then it contains an edge $e=[a,c,b]$ or $e=[b,c,a]$, where $a$ is a boundary vertex and $a$ does not appear as an edge label in $\cG$. We remove the interior of $e$ and the vertex $a$ from $\cG$ to obtain a labeled oriented tree $\cG'$ with fewer vertices. As before,  
%on the presentation level 
this amounts 
%to a $Q^{**}$-transformation and hence 
to a 3-deformation $h\colon K(\cG) \to K(\cG')$ on the 2-complex level. If $\cG$ is interior reducible at a vertex $b$ then there exist edges $e_1=[a,d,b]$ and $e_2=[c,d,b]$, or $e_1=[b,d,a]$ and $e_2=[b,d,c]$. Fold the two edges $e_1$ and $e_2$ into one edge and label it by $d$, matching orientations and identifying the vertices $a$ and $c$ into one vertex $a$. Change all edge labels $c$ into edge labels $a$. Again,  
%on the presentation level 
this amounts 
%to a $Q^{**}$-transformation and hence 
to a 3-deformation on the 2-complex level. 
  
If $\cH$ is a sub-LOT of $\cG$ then the reductions just described transform $\cH$ into a sub-LOT $\cH_0$ of $\cG_{red}$ and the homotopy equivalence $h\colon K(\cG) \to K(\cG_{red})$ restricts to a homotopy equivalence $h\colon K(\cH) \to K(\cH_0)$.

In summary we have the following result that will be used later in this article.

\lemanf\label{LOTtransformation} Let $\cal G$ be a labeled oriented tree and $\cal H$ be a sub-LOT. Transform $\cal G$ to ${\cal G}_{red}$ and let ${\cal H}_0$ be the image of $\cal H$ in ${\cal G}_{red}$ under that transformation. Then there is a commutative diagram

\begin{center}
\begin{tikzpicture}[scale = 0.7]
\node [ ] at (0,2) {$K({\cal G})$};
\node [ ] at (3,2) {$K({\cal G}_{red})$};
\node [ ] at (0,0) {$K({\cal H})$};
\node [ ] at (3,0) {$K({\cal H}_0)$};
\draw [->] (1,2) to (2,2);
\draw [->] (1,0) to (2,0);
\draw [<-] (0,1.5) to (0,0.5);
\draw [<-] (3,1.5) to (3,0.5);
\node [above] at (1.5,2) {$h$};
\node [left] at (0,1) {$i$};
\node [right] at (3,1) {$i$};
\node [above] at (1.5,0) {$h_0$};
\end{tikzpicture}
\end{center}

where the horizontal maps are homotopy equivalences and the vertical maps are inclusions. 
\lemende

We conclude this section with an outline of the paper. In Section \ref{sec:combinatorialasphericity} we review some basic concepts from combinatorial topology and introduce the notion of relative vertex asphericity. Section \ref{sec:relative} contains a test for relative vertex asphericity based on a result of Stallings. The material in this section is of independent interest with possible applications not directly connected with the study of labeled oriented trees.  Section \ref{sec:alter} introduces altered LOT-presentations and contains Theorem \ref{sprgr}, a result shown by Huck and the second author \cite{HR01}. There it was used to show that prime injective labeled oriented trees are aspherical (prime means that the labeled oriented tree does not contain sub-LOTs). Corollary \ref{orientationcorollary} is a relative version of Theorem \ref{sprgr}. Section \ref{mainproof} contains the proof of our main result Theorem \ref{hauptsatz}. We also show that if ${\cal H}$ is 
 a sub-LOT of an injective labeled oriented tree $\cG$, then $G({\cal H})$ is a subgroup of $G(\cG)$. In Section \ref{sec:consequences} we extend some of our results to a class of non-injective labeled oriented trees, and we close with an application to virtual knots.

\section{Relative vertex asphericity}\label{sec:combinatorialasphericity}

In the following we work in the category of combinatorial 2-complexes and combinatorial maps. Recall that a map $f\colon X\to Y$ between CW complexes is said to be {\em combinatorial} if the restriction of $f$ to each open cell of $X$ is a homeomorphism onto its image. Most of the 2-complexes considered in this article will be 
standard 2-complexes built from
group presentations. Such 2-complexes have a single $0$-cell.  

\defanf\label{sphericaldiag1} A {\em spherical diagram} over a 2-complex $K$ is a combinatorial map \map{C}{K}, where $C$ is a cell decomposition of
the 2-sphere.
\defende

If the 2-complex $K$ is the standard 2-complex $K(P)$ associated with a group presentation $P=\langle {\bf x}\ |\ {\bf r} \rangle$ one can define spherical diagrams in graph theoretic terms. Note that $K(P)$ has a single 0-cell and its oriented 1-cells are in one-to-one correspondence with elements from ${\bf x}$. Thus there is a one-to-one correspondence between edge paths in $K(P)$ and words in ${\bf x}^{\pm 1}$. Suppose $f\colon C\to K(P)$ is a spherical diagram. If $e$ is an edge in $C$ that gets mapped to the edge $x\in {\bf x}$ in $K(P)$, then we orient $e$ to make $f$ orientation preserving on $e$ and label it by $x$. In this way the spherical diagram gives rise to a planar connected oriented graph $C^{(1)}$ such that
\begin{itemize}
\item oriented edges in $C^{(1)}$ are labeled by elements from ${\bf x}$;
\item the word read off the boundary path of an inner or the outer region of $C^{(1)}$ is a cyclic permutation of a word $r^{\epsilon}$, where $\epsilon=\pm 1$ and $r\in{\bf r}$.
\end{itemize}
We can also start with a planar connected oriented graph $C^{(1)}$, labeled in the fashion just described, and produce a spherical diagram $f\colon C\to K(P)$. This provides an alternate definition for spherical diagrams over standard 2-complexes $K(P)$ in terms of planar connected oriented graphs with labeled edges.
This purely combinatorial point of view is taken by Bogley and Pride \cite{BogleyPride}. They work with pictures over group presentations which are obtained by dualizing spherical diagrams. Another standard reference for diagrams is Gersten \cite{Ger87}. 

\defanf Let $\Gamma$ be a graph and $\Gamma_0$ be a subgraph (which could be empty).
\begin{enumerate}
\item 
An edge cycle $c=e_1\ldots e_q$ in $\Gamma$ is called {\em homology reduced} if it does not contain a pair of edges $e_i$ and $e_j$ so that $e_j=\bar e_i$, where $\bar e_i$ is the edge $e_i$ with opposite orientation.
\item An edge cycle $c=e_1\ldots e_q$ is said to be {\em homology reduced relative to $\Gamma_0$} if it does not contain a pair of edges $e_i$ and $e_j$ of $\Gamma - \Gamma_0$ so that $e_j=\bar e_i$.
\end{enumerate}
\defende 

If $v$ is a vertex of a 2-complex $K$ then the link $Lk(K,v)$ is the boundary of a regular neighborhood of $v$ in $K$ equipped with the induced cell decomposition. So $Lk(K,v)$ is a graph. If $K$ has a single vertex we denote by $Lk(K)$ the link of that vertex. We refer to the edges in $Lk(K)$ as {\em corners}, since they can be thought of as the corners of the 2-cells of $K$.\\

Let $K$ be a 2-complex with a single vertex and $K_0$ a subcomplex (which could be empty).  Note that $Lk(K_0)$ is a subgraph of $Lk(K)$. Let $f\colon C\to K$ be a spherical diagram and $v$ be a vertex of the 2-sphere $C$. The map $f$ induces a combinatorial map $f_L\colon Lk(C,v)\to Lk(K)$. Note that $Lk(C,v)$ is a circle and the image of that circle, oriented clockwise, is a cycle of corners $\alpha(v)=\alpha_1\ldots \alpha_q$, that is a closed edge path, in $Lk(K)$. We say that the diagram $f\colon C\to K$ is {\em vertex reduced at $v$ relative to $K_0$} if the cycle $\alpha(v)$ in $Lk(K)$ is homology reduced relative to $Lk(K_0)$. 
We say that the diagram is {\em vertex reduced relative to $K_0$} if it is vertex reduced relative to $K_0$ at all its vertices. The 2-complex $K$ is called {\em vertex aspherical relative to $K_0$}, VA relative to $K_0$ for short, if for every spherical diagram $f\colon C\to K$ that is vertex reduced relative to $K_0$ we have $f(C)\subseteq K_0$. 

If we omit  ``relative to" we implicitly imply relative to the empty set $\emptyset$, even if a subcomplex is present. For example if we say a spherical diagram is vertex reduced we mean vertex reduced relative to $\emptyset$.

\pagebreak

\satzanf\label{sVA} If $K$ is VA relative to $K_0$, then $\pi_2(K)$ is generated, as $\pi_1(K)$-module, by the image of $\pi_2(K_0)$ under the map induced by inclusion. In particular, if $K_0$ is aspherical, then so is $K$.
\satzende

\bewanf This follows from the fact that $\pi_2(K)$ is generated by spherical
diagrams (see \cite{BogleyPride} Theorem 1.3, p.\ 162 and the literature cited there) and the fact that each spherical diagram $f\colon C\to K$ is homotopic  to a spherical diagram $f'\colon C'\to K$ that is vertex reduced.
\bewende

Vertex asphericity in case $K_0=\emptyset$ was considered in Huck, Rosebrock \cite{HR02}. 
It presents a generalization of {\em diagrammatic reducibility}, DR for short. See Gersten \cite{Ger87} for a definition. We note that diagrammatic reducibility implies vertex asphericity, and vertex asphericity implies asphericity (see \cite{HR02}). \\

\section{\label{sec:relative} A test for relative vertex asphericity}

A graph is called a {\em forest} if its connected components are trees.

\defanf Let $\Gamma$ be a graph and $\Gamma_0$ be a subgraph. 
\begin{itemize} 
\item The graph $\Gamma$ is called a {\em forest relative to $\Gamma_0$} if every homology reduced cycle is contained in $\Gamma_0$.
\item The graph $\Gamma$ is called a {\em tree relative to $\Gamma_0$} if $\Gamma $ is connected and every homology reduced cycle is contained in $\Gamma_0$.
\end{itemize}
\defende

\lemanf \label{graphlemma} Let $\Gamma$ be a graph, $\Gamma_0$ a subgraph with connected components $\Gamma_1,\ldots ,\Gamma_n$. Let $\Gamma'$ be the graph obtained by collapsing each component $\Gamma_i$ to a vertex $g_i\in \Gamma_i$. Then $\Gamma$ is a forest relative to $\Gamma_0$ if and only if $\Gamma'$ is a forest.
\lemende

\bewanf If $\Gamma $ is a forest relative to $\Gamma_0$
then each homology reduced cycle is contained in some $\Gamma_i$. After collapsing
each $\Gamma_i$ to a point there will be no more non-trivial homology reduced
cycles and so $\Gamma'$ is a forest.

For the converse let $Y$\ be the closure of $\Gamma - \Gamma_0$ in $\Gamma$. The intersection $\Gamma_0\cap Y$ is a set of vertices. 
Observe that since $\Gamma'$ is a forest, if an edge $e$ of $Y$ appears in a cycle in $\Gamma$, then so must its inverse $\bar e$. 
In particular, a homology reduced cycle in $\Gamma$ has to be
contained in $\Gamma_0$ and hence in one of the connected components $\Gamma_i$.
\bewende

We next discuss a graph theoretic result. 
Let $C$ be a cell decomposition of an oriented  2-sphere with oriented edges. A {\em sink} is a vertex in $C$ with all adjacent edges pointing towards it, a {\em source} is a vertex in $C$ with all adjacent edges pointing away from it. A 2-cell is {\em consistently oriented} if all its boundary edges are oriented clockwise or all are oriented anti-clockwise. We say a 2-cell has {\em exponent sum zero} if, when reading around its boundary in clockwise direction, one encounters the same number of clockwise oriented edges as anti-clockwise oriented edges. The lemma below is due to Stallings \cite{Sta87},  Lemma 1.2.

\lemanf\label{lstallings} Given a cell decomposition $C$ of the 2-sphere with oriented edges.
If $C$ does not contain a consistently oriented 2-cell, then it contains a sink or a source.
\lemende

Let $P=\langle {\bf x}\ |\ {\bf r}\rangle$ be a presentation. The link of the single vertex in the associated standard 2-complex is also referred to as the {\em Whitehead graph} $W(P)$.
It can be defined directly from the presentation without reference to topological notions. See for example \cite{BogleyPride}, page 170. The Whitehead graph $W(P)$ is a non-oriented graph on vertices $\{ x^+, x^- \ |\ x\in {\bf x}\}$, where $x^+$ is a point of the oriented edge $x$ of $K(P)$ close to the beginning of that edge, and $x^-$ is a point close to the ending of that edge. Vertices $x^{\epsilon}$ and $y^{\delta}$, ($x,y\in{\bf x},\,\, \epsilon, \delta \in \{ \pm  \}$),  are connected by an edge in $W(P)$ if there is a 2-cell in $K(P)$ with a corner connecting the two points. For that reason we refer to the edges of $W(P)$ also as {\em corners}. The {\em positive graph} $W^+(P)\subset W(P)$  is the
full subgraph on the vertex set  $\{ x^+ \ |\ x\in {\bf x}\}$, the {\em
negative graph } $W^-(P)\subset W(P)$ is the full subgraph on the
vertex set $\{x^- \ |\ x\in {\bf x}\}$.\\

\defanf\label{Stdef} A presentation $P$ is said to {\em satisfy the Stallings test} if the following conditions hold:
\begin{itemize} 
\item Relator conditions:
\begin{enumerate}
\item[(a)] Relators of $P$ are cyclically reduced.
\item[(b)] Relators of $P$ are not positive or negative words.
\end{enumerate}
\item Forest  condition:
\begin{enumerate}
\item [(c)]$W^+(P)$ and $W^-(P)$  are forests.
\end{enumerate}
\end{itemize}
\defende

The following application of Lemma \ref{lstallings} is well known. It first appeared in
Gersten \cite{Ger87}, Proposition 4.12, in the context of Adian presentations.

\satzanf\label{StallingsTest} Let $P$ be a presentation that satisfies the Stallings test. Then $K(P)$ is VA (even DR) and hence aspherical.
\satzende

The main result of this section is a relative version of Theorem \ref{StallingsTest}. Before we can state it we need some notation. If $P_1=\langle {\bf x}_1\ |\ {\bf r}_1\rangle$ and $P_2=\langle {\bf x}_2\ |\ {\bf r}_2\rangle$ are presentations then $P_1\cup P_2=\langle {\bf x}_1\cup {\bf x}_2\ |\ {\bf r}_1\cup {\bf r}_2\rangle$. Let $P=\langle {\bf x}\ |\ {\bf r} \rangle$ be a presentation and let $\{ T_1, \ldots ,T_n\}$ be a set of disjoint {\em full} sub-presentations of $P$. Full means that if $r$ is a relator in $P$ that only involves generators from $T_i$, then $r$ is already a relator in $T_i$.  Disjoint means that the generating sets 
of $T_i$ and $T_j$ are disjoint subsets of ${\bf x}$ in case $i\ne j$. Let $T=T_1\cup \ldots \cup T_n$. The complex $K(T)=K(T_1)\vee \ldots \vee K(T_n)$  is a sub-complex of $K(P)$.

Let $T_i=\langle {\bf t}_i\ |\ {\bf s}_i \rangle$ and let ${\bf u}_i$ be the set of words with letters in ${\bf t}_i^{\pm 1}$ of exponent sum zero. Note that we do not assume that the words in ${\bf u}_i$ are freely or cyclically reduced.  Let ${\bf T}_i=\langle {\bf t}_i\ |\ {\bf s}_i\cup{\bf u}_i \rangle$. Let ${\bf T}={\bf T}_1\cup \ldots \cup {\bf T}_n$ and note that $P\cup {\bf T}=\langle {\bf x} \ |\ {\bf r}\cup{\bf u}_1\cup \ldots \cup{\bf u}_n \rangle$ since ${\bf t}_i\subseteq {\bf x}$ and ${\bf s}_i\subseteq {\bf r}$ for every $i$. The presentation $P\cup{\bf T}$ is infinite and in the group $G(P\cup{\bf T})$ the generators of each $T_i$ are identified, since we have the relator $t^{-1}t'$ in $P\cup{\bf T}$ for every pair $t, t'$ of generators in $T_i$.  Note that ${\bf T}_i$ is a sub-presentation of $P\cup{\bf T}$ and the subgraph $W({\bf T}_i)$ of the Whitehead graph $W(P\cup{\bf T})$, which is spanned by the vertices $t^{\pm}$ with $t\in {\bf t}_i$, contains the complete graph on these vertices.  In fact, every pair of vertices in $W({\bf T}_i)$ is connected by infinitely many edges, and at every vertex in $W({\bf T}_i)$ there are attached infinitely many loops. 

\defanf\label{relStdef} Let $P$ be a presentation, and 
let $\{ T_1, \ldots ,T_n\}$ be a set of disjoint full sub-presentations. Let $T=T_1\cup \ldots \cup T_n$. Then $P$ is said to {\em satisfy the Stallings test relative to $T$} if the following conditions hold:
\begin{itemize}
\item Relator conditions:
\begin{enumerate}
\item[(a)] Relators of $P-T$ are cyclically reduced.
\item[(b)] Relators of $P-T$ are not positive or negative words.
\item[(c)] Relators of $T$ have exponent sum zero.
\item[(d)] Any word $w$ in the generators of some $T_i$ that represents the trivial element of the group defined by $P\cup {\bf T}$ has exponent sum zero.
\end{enumerate}
\end{itemize}
\begin{itemize}
\item Forest condition:
\begin{enumerate}
\item[(e)] $W^+(P\cup {\bf T})$ is a forest relative to $W^+({\bf T})$ and $W^-(P\cup {\bf T})$ is a forest relative to $W^-({\bf T})$.
\end{enumerate}
\end{itemize}
\defende

Here is the main result of this section:

\satzanf\label{relativeStallingstest} If $P$ is a presentation that satisfies the Stallings test relative to $T$, then $K(P)$ is VA relative to $K(T)$. Furthermore, the inclusion induced homomorphism $\pi_1(K(T_i))\to \pi_1(K(P))$ is injective for every $i=1,\ldots ,n$.
\satzende

For the proof of this theorem we need the following lemma:

\lemanf\label{lVA}
Let $P$ be a presentation and 
let $\{ T_1, \ldots ,T_n\}$ be a set of disjoint full sub-presentations, $T=T_1\cup \ldots \cup T_n$. Assume that the relators of $P$ satisfy the relator conditions (a)-(d) of the relative Stallings test (Definition \ref{relStdef}). 
If $K(P)$ is not VA relative to $K(T)$, then there is a spherical diagram $f\colon C\to K(P\cup {\bf T})$ such that
\begin{enumerate}
\item $f\colon C\to K(P\cup {\bf T})$ is vertex reduced relative to $K({\bf T})$,
\item $f(C)$ is not contained in $K({\bf T})$, and
\item If $v$ is a vertex in $C$ then the corner cycle $\alpha(v)$ has length at least two and does not contain two distinct corners $\alpha_p$ and $\alpha_q$ that both come from relators of one 
subpresentation ${\bf T}_i$ of $P\cup {\bf T}$, $i\in\{1,\ldots ,n\}$.
\end{enumerate}
\lemende

\bewanf Since  $K(P)$ is not VA relative to $K(T)$ there exists a spherical diagram $C^*\to K(P)$ that is vertex reduced relative to $K(T)$ but does not map entirely into $K(T)$.
Since $K(P)$ is a subcomplex of $K(P\cup {\bf T})$ this diagram can be viewed as a diagram $C^* \to K(P\cup {\bf T})$ 
which is vertex reduced relative to $K({\bf T})$ but does not map entirely into $K({\bf T})$. Let $\Omega$ be the collection of all spherical diagrams that have that feature. Consider the subset $\Omega_0\subseteq \Omega$ of those for which the 2-sphere $C$ contains the smallest number of 2-cells. From $\Omega_0$ choose a spherical diagram $f\colon C\to K(P\cup {\bf T})$ for which $C$ has the smallest number of edges.

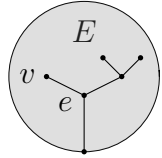
\begin{figure}[ht]\centering
\begin{tikzpicture}[scale=0.5]

%% 2-cell
\draw (2,0) arc (0:365:2);
%%% 2-cell label
\node at (0,1.2){$E$};
%%% 2-cell color
\fill[gray,nearly transparent] (2,0) arc (0:365:2);

%% vertices
\fill (0,-2) circle (2pt);
\fill (0,-0.5) circle (2pt);
\fill (-1,0) circle (2pt); \node[left] at (-1,0){$v$};
\fill (1,0) circle (2pt);
\fill (0.5,0.5) circle (2pt);
\fill (1.5,0.5) circle (2pt);

%% 1-cells
\draw (0,-2)--(0,-0.5);
\draw (0,-0.5)--(-1,0); \node[below] at (-0.5,-0.25){$e$};
\draw (0,-0.5)--(1,0);
\draw (1,0)--(0.5,0.5);
\draw (1,0)--(1.5,0.5);

\end{tikzpicture}

\caption{\label{valency1} Vertices of valency one in $C$.}
\end{figure}

This spherical diagram does not contain a vertex of valency one. If $v$ were a vertex of valency one in $C$, then it would be a vertex in the boundary of a cell $E$ that maps to some $K({\bf T}_i)$ since we assumed the relators in $P$ outside $T$ to be cyclically reduced. Let $e$ be the edge in $C$ that contains $v$. See Figure \ref{valency1}. We can remove $v$ and the interior of $e$ and transform $E$ into $E'$. Note that the boundary words of $E$ and $E'$ are the same up to free or cyclic reduction, hence removing $v$ and the interior of $e$ produces a spherical diagram $f'\colon C'\to K(P\cup {\bf T})$ with fewer edges but the same number of 2-cells contradicting the choice of $f$.
%%%%%%%%%%%%%%%%%%%%%%%%%%%%
% tree in disk
%%%%%%%%%%%%%%%%%%%%%%%%%%%%
%%%%%%%%%%%%%%%%%%%%%%%%%%%%%%%%%%%%%%%
%%%%%%%%%%%%%%%%%%%%%%%%%%%%%%%%%%%%%%%

\begin{figure}[ht]\centering
\begin{tikzpicture}

% left diagram

%% 2-cells
\draw (-3,0)--(-4,1)--(-5,0)--(-4,-1)--cycle;
\draw (-3,0)--(-2,1)--(-1,0)--(-2,-1)--cycle;
\draw (-3,0)--(-cos 60 - 3,1); \draw (-3,0)--(cos 60 -3,1); \draw (-3,0)--(-3,-1);
%%% 2-cell labels
\node at (-4,0){$E_p$};
\node at (-2,0){$E_q$};
%%% 2-cell colors
\fill[yellow,nearly transparent] (-3,0)--(-4,1)--(-5,0)--(-4,-1)--cycle;
\fill[cyan,nearly transparent] (-3,0)--(-2,1)--(-1,0)--(-2,-1)--cycle;

%% vertex
\fill (-3,0) circle (2pt);

%% links \& labels
\draw (-3 -0.6*cos 45, 0.6*sin 45) arc (135:225:0.6);
\node[left] at (-3,0){$\alpha_p$};
\draw (-3 + 0.6*cos 45, 0.6*sin 45) arc (45:-45:0.6);
\node[right] at (-3,0){$\alpha_q$};

%right diagram

%% 2-cell
\draw (3,0.3)--(4,1)--(5,0)--(4,-1)--(3,-0.3)--(2,-1)--(1,0)--(2,1)--cycle;
\draw (3,0.3)--(-cos 60 + 3,1); \draw (3,0.3)--(cos 60 + 3,1); \draw (3,-0.3)--(3,-1);
%%% 2-cell label
\node at (3,0.01){$E$};
%% 2-cell colors
%\fill[cyan,nearly transparent] (247,0)--(3,0.3)--(4,1)--(5,0)--(4,-1)--(3,-0.3)--cycle;
%\fill[yellow,nearly transparent] (187,0)--(3,0.3)--(2,1)--(1,0)--(2,-1)--(3,-0.3)--cycle;
%\fill[green,nearly transparent] (187,0)--(3,0.3)--(247,0)--(3,-0.3)--cycle;
\fill[green,nearly transparent] (3,0.3)--(4,1)--(5,0)--(4,-1)--(3,-0.3)--(2,-1)--(1,0)--(2,1)--cycle;

%% vertices
\fill (3,0.3) circle (2pt); \fill (3,-0.3) circle (2pt);

\end{tikzpicture}

\caption{\label{split} Splitting at a vertex.}
\end{figure}
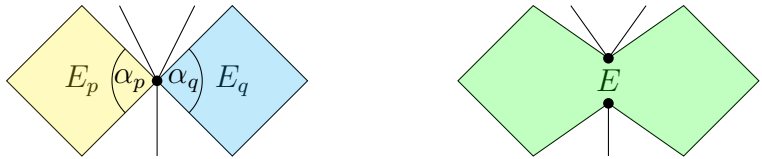

The spherical diagram under consideration has the first and the second property by choice. Let us look at the third property. Let $v\in C$ be a vertex and $\alpha(v)$ its corner cycle.
The length of $\alpha(v)$ is at least two because $C$ does not contain vertices of valency one. Suppose $E_p, E_q$ are 2-cells in $C$ that both get mapped to 2-cells of $K({\bf T}_i)$ and give rise to the corners $\alpha_p, \alpha_q\in\alpha(v)$, respectively. Let us assume first that $E_p$ and $E_q$ are distinct (we do not rule out that $E_p$ and $E_q$ share boundary edges). We can split $v$ in $C$ into two vertices, and fuse the cells $E_p$ and $E_q$ into one cell $E$. See Figure \ref{split}. 
%%%%%%%%%%%%%%%%%%%%%%%%%%%%
% merging squares
%%%%%%%%%%%%%%%%%%%%%%%%%%%%
%%%%%%%%%%%%%%%%%%%%%%%%%%%%%%%%%%%%%%%
%%%%%%%%%%%%%%%%%%%%%%%%%%%%%%%%%%%%%%%
This creates a new cell division $C'$ with one fewer 2-cell than $C$. Note that the boundary word of $E$ is a word obtained by reading around the boundary of the union $E_p\cup E_q$, thus it is a word in ${\bf t}_i^{\pm 1}$, where ${\bf t}_i$ is the generating set of $T_i$, of exponent sum zero. So the boundary word of $E$ is a relator of ${\bf T}_i$. The splitting process has created a new spherical diagram $f'\colon C'\to K(P\cup {\bf T})$. It is vertex reduced relative to $K({\bf T})$ and $f'(C')$ is not contained in $K({\bf T})$. Since $C'$ contains fewer 2-cells than $C$ this contradicts the choice of $f\colon C\to K(P\cup {\bf T})$.\\

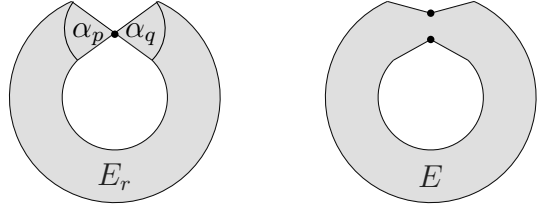
\begin{figure}[ht]\centering

\begin{tikzpicture}[scale=0.7]

% left diagram

%% 2-cell
\draw (-cos 45 - 3, sin 45) arc (135:405:1);
\draw (-2*cos 65.705 -3, 2*sin 65.705) arc (114.295:425.705:2);
\draw (-cos 45 - 3, sin 45)--(2*cos 65.705 - 3, 2*sin 65.705);
\draw (cos 45 - 3, sin 45)--(-2*cos 65.705 - 3, 2*sin 65.705);
%%% 2-cell label
\node at (-3,-1.5){$E_r$};
%%% 2-cell color
\fill[gray,nearly transparent] (-3,1.2)--(-2*cos 65.705 -3, 2*sin 65.705) arc (114.295:425.705:2)--(2*cos 65.705 - 3, 2*sin 65.705)--(-3,1.2)--(cos 45 - 3, sin 45) arc (45:-225:1)--(-cos 45 - 3, sin 45)--cycle;

%% vertex
\fill (-3,1.2) circle (2pt);

%% arcs & lables
\draw (cos 45 -3, sin 45) arc (-45:37:.85);
\node at (-2.5,1.2){$\alpha_q$};
\draw (-cos 45 - 3, sin 45) arc (225:143:.85);
\node at (-3.5,1.2){$\alpha_p$};

% right diagram

%% 2-cell
\draw (-cos 45 + 3, sin 45) arc (135:405:1);
\draw (-2*cos 65.705 + 3, 2*sin 65.705) arc (114.295:425.705:2);
\draw (-cos 45 + 3, sin 45)--(3, 1.1);
\draw (cos 45 + 3, sin 45)--(3, 1.1);
\draw (2*cos 65.705 + 3, 2*sin 65.705)--(3, 1.6);
\draw (-2*cos 65.705 + 3, 2*sin 65.705)--(3, 1.6);
%%% 2-cell label
\node at (3,-1.5){$E$};
%%% 2-cell color
\fill[gray,nearly transparent] (3,1.6)--(-2*cos 65.705 + 3, 2*sin 65.705) arc (114.295:425.705:2)--(2*cos 65.705 + 3, 2*sin 65.705)--(3,1.6)--(3,1.1)--(cos 45 + 3, sin 45) arc (45:-225:1)--(-cos 45 + 3, sin 45)--(3,1.1)--cycle;

%% vertices
\fill (3,1.1) circle (2pt);
\fill (3,1.6) circle (2pt);

\end{tikzpicture}
\caption{\label{split2} Splitting at a vertex.}
\end{figure}

We next assume that there is one 2-cell $E_r$ at $v$ in $C$ that maps to a 2-cell in ${\bf T}_i$ that gives rise to both corners $\alpha_p$ and $\alpha_q$. See Figure \ref{split2}. 
%%%%%%%%%%%%%%%%%%%%%%%%%%%%
% annulus
%%%%%%%%%%%%%%%%%%%%%%%%%%%%
%%%%%%%%%%%%%%%%%%%%%%%%%%%%%%%%%%%%%%%
%%%%%%%%%%%%%%%%%%%%%%%%%%%%%%%%%%%%%%%
We can split at $v$ as before, but note that now $C'$ is not a cell division of the 2-sphere because
the 1-skeleton $C'^{(1)}$ is not connected. Let $C'^{(1)}_1$ and $C'^{(1)}_2$ be the connected components of $C'^{(1)}$. The boundary word for the outer region of each planar connected graph $C'^{(1)}_j$, $j=1,2$, is a word in ${\bf t_i}^{\pm 1}$ which represents the trivial element of the group defined by $P\cup {\bf T}$ and hence (relator condition (d) of Definition \ref{relStdef}) has exponent sum zero. Thus the boundary word is a relator in ${\bf T}_i$. It follows that each $C'^{(1)}_j$ determines a spherical diagram $f_j\colon C'_j\to K(P\cup {\bf T})$ that is vertex reduced relative to $K({\bf T})$. Since $f(C)$ is not in $K({\bf T})$, the image $f_j(C'_j)$ is not in $K({\bf T})$ for one of the $j$'s. But $C'_j$ has fewer 2-cells than $C$. This contradicts the choice of $f\colon C\to K(P\cup {\bf T})$.

This shows that the spherical diagram $f\colon C\to K(P\cup {\bf T})$ has the three properties stated above.
\bewende

{\bf Proof of Theorem \ref{relativeStallingstest}}: Suppose $K(P)$ is not VA relative to $K(T)$. Then there exists a spherical diagram $f\colon C\to K(P\cup {\bf T})$ that satisfies the conditions 1, 2, and 3 stated in Lemma \ref{lVA}. The relator conditions (b) and (c) in Definition \ref{relStdef} (relative Stallings test) imply that $C$ does not contain cells with consistently oriented boundary, hence $C$ contains a sink or a source by Lemma \ref{lstallings}. Let us assume without loss of generality that $C$ contains a source, say at the vertex $v\in C$. The cycle $\alpha(v)=\alpha_1\ldots \alpha_l$
satisfies $l\ge 2$, is contained in $W^+(P\cup {\bf T})$ and is homology reduced relative to $W^+({\bf T})$ because $f\colon C \to K(P\cup {\bf T})$ is vertex reduced relative to $K({\bf T})$. Since $W^+(P\cup {\bf T})$ is a forest relative to $W^+({\bf T})$ we know that $\alpha(v)$
is entirely contained in a connected component of $W^+({\bf T})$, and hence in some $W^+({\bf T}_i)$, because $W^+({\bf T})$ is a disjoint union of the $W^+({\bf T}_i)$, $i=1,\ldots ,n$. Thus, if $\alpha(v)=\alpha_1\ldots\alpha_l$, then all corners $\alpha_j$, $j=1,\ldots ,l$ are in $W^+({\bf T}_i)$. This contradicts condition 3 of Lemma \ref{lVA}. We have reached a contradiction. Thus $K(P)$ is VA relative to $K(T)$.

Suppose the map $\pi_1(K(T_i))\to \pi_1(K(P))$ is not injective for some $i=1,\ldots ,n$. Then there exists a cyclically reduced 
word $w$ in the generators of $T_i$ that represents the trivial element of $\pi_1(K(P))$ but is not trivial in 
$\pi_1(K(T))$. Hence $w$ is the boundary word of a vertex reduced disc diagram $g\colon D\to K(P)$. Note that $D$ has to contain 2-cells that are not mapped to $K(T)$ because the map $\pi_1(K(T_i))\to \pi_1(K(T))=\pi_1(K(T_1))*\ldots *\pi_1(K(T_n))$ is injective.  By the relator condition (d) of Definition \ref{relStdef} (relative Stallings test) the word $w$ has exponent sum zero and hence is a relator of ${\bf T}_i$.  We can attach a disc $D'$ to $D$ and obtain a spherical diagram $f\colon C\to K(P\cup {\bf T})$. Note that this spherical diagram is vertex reduced.  If it were not, then there would have to be a vertex on the boundary of $D$ where the spherical diagram $f\colon C\to K(P\cup {\bf T})$ is not vertex reduced. But that would mean that $D$ contains a 2-cell with boundary word $w$. This would imply that $w$ is a relator in $P$. Since we assumed $T_i$ to be a full sub-presentation, the word $w$ would have to be a relator in $T_i$, which is not the case because $w$ does not represent the trivial element of $\pi_1(K(T_i))$. Thus $f\colon C\to K(P\cup {\bf T})$ is indeed vertex reduced. Since $K(P)$ is VA relative to $K(T)$ we have that $f(C)\subseteq K({\bf T})$, which implies that $g(D)\subseteq K(P)\cap K({\bf T})=K(T)$. We have reached a contradiction.\bewende

\section{\label{sec:alter}Altering LOT-presentations and orienta\-tions} 

We say a labeled oriented graph $\cal Q$ is a {\em reorientation} of a labeled oriented graph $\cal P$ if $\cal Q$ is obtained from $\cal P$ by changing the orientation of each edge of a subset of the set of edges of $\cal P$.\\

The next result is proved in Section 3 of \cite{HR01}. 

\satzanf\label{sprgr}  Let $\cal P$ be a compressed injective labeled oriented tree
that does not contain a boundary reducible sub-LOT.  Then there is a
reorientation $\cal Q$ of $\cal P$ such that $W^+(Q)$ and $W^-(Q)$ are trees, where $Q$ is the LOT-presentation associated with ${\cal Q}$.
\satzende

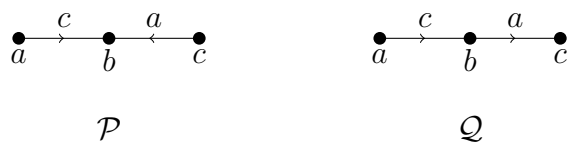
\begin{figure}[ht]\centering
\begin{tikzpicture}[scale=1.2]

\fill (0,0) circle (2pt);
\node[below] at (0,0) {$a$};
\fill (1,0) circle (2pt);
\node[below] at (1,0) {$b$};
\fill (2,0) circle (2pt);
\node[below] at (2,0) {$c$};

\begin{scope}[decoration={markings, mark=at position 0.5 with {\arrow{>}}}]
\draw [postaction={decorate}] (0,0) -- (1,0) node[midway, above]{$c$};
\end{scope}

\begin{scope}[decoration={markings, mark=at position 0.5 with {\arrow{<}}}]
\draw [postaction={decorate}] (1,0) -- (2,0) node[midway, above]{$a$};
\end{scope}

\node at (1,-1) {$\mathcal P$};

\fill (4,0) circle (2pt);
\node[below] at (4,0) {$a$};
\fill (5,0) circle (2pt);
\node[below] at (5,0) {$b$};
\fill (6,0) circle (2pt);
\node[below] at (6,0) {$c$};

\begin{scope}[decoration={markings, mark=at position 0.5 with {\arrow{>}}}]
\draw [postaction={decorate}] (4,0) -- (5,0) node[midway, above]{$c$};
\end{scope}

\begin{scope}[decoration={markings, mark=at position 0.5 with {\arrow{>}}}]
\draw [postaction={decorate}] (5,0) -- (6,0) node[midway, above]{$a$};
\end{scope}

\node at (5,-1) {$\mathcal Q$};

\end{tikzpicture}
\caption{\label{reorientexample} A reorientation $\cQ$ of a labeled oriented tree $\cP$.}
\end{figure}

\bspanf\label{firstexample} Consider the compressed injective labeled oriented tree $\cP$ shown in Figure \ref{reorientexample}. Let $P$ be its LOT-presentation. Then $W^-(P)$ is a tree but $W^+(P)$ is not. However, for the reorientation $\cQ$ of $\cP$ (also shown in Figure \ref{reorientexample}) both $W^-(Q)$ and $W^+(Q)$ are trees, where $Q$ is the LOT-presentation of $\cQ$.
\bspende

The situation is more complicated in the presence of a boundary reducible sub-LOT. In fact, if $\cP$ is a labeled oriented tree and $\cT$ is a sub-LOT that is not boundary reduced, then there does not exist a reorientation $\cQ$ of $\cP$ so that both $W^-(Q)$ and $W^+(Q)$ are forests. This follows from the following lemma. 

\lemanf\label{notboundaryredlemma} If $\cT$ is a labeled oriented tree that is not boundary reduced, then either $W^-(T)$ or $W^+(T)$ is not a forest, where $T$ is the LOT-presentation of $\cT$.
\lemende

\bewanf Since $\cT$ is not boundary reduced there exists a boundary vertex $a$ that does not occur as edge label. Thus, depending on the orientation of the edge of $\cT$ containing $a$, either $a^+$ is an isolated vertex in $W^+(T)$, or $a^-$ is an isolated vertex in $W^-(T)$ (a vertex in a graph is {\em isolated} if it is not the vertex of an edge). Let us assume without loss of generality that $a^+$ is isolated in $W^+(T)$. If $n$ is the number of vertices in $\cT$, then $W^+(T)$ contains $n$ vertices and $n-1$ edges. Since $a^+$ is isolated, $W^+(T)$ contains a subgraph containing $n-1$ vertices and $n-1$ edges. Such a graph contains a cycle.
\bewende

If $\cP$ is a labeled oriented tree that contains a sub-LOT $\cT$ that is not boundary reduced, then a reorientation $\cQ$
of $\cP$ contains a reorientation $\cT^*$ of  $\cT$, which is also not boundary reduced. By Lemma \ref{notboundaryredlemma} either $W^-(\cT^*)$ or $W^+(\cT^*)$ is not a tree. Since $W^-(\cT^*)\subseteq W^-(\cQ)$ and $W^+(\cT^*)\subseteq W^+(\cQ)$, it follows that either $W^-(\cQ)$ or $W^+(\cQ)$ is not a tree.
Thus, Theorem \ref{sprgr} does not hold in the presence of a boundary reducible sub-LOT. However, a relative version of that theorem does hold. Before we state it, we make some assumptions and definitions.\\

We assume for the remainder of this section: $\cP$ is a reduced labeled oriented tree and $\{{\mathcal T}_1,\ldots ,{\mathcal T}_n \}$ is a set of proper maximal sub-LOTs.  The sub-LOTs are pairwise disjoint, that is  $\cT_i\cap\cT_j=\emptyset$ in case $i\ne j$. Let  $\cT=\cT_1\cup...\cup\cT_n$. Let $P$, $T$, $T_i$ be the LOT-presentations associated with $\cP$, $\cT$, $\cT_i$, respectively.\\

We denote by $\cP-\cT$ the forest with edge set the edges in $\cP$ not in $\cT$
and with vertex set the vertices which bound the edges in $\cP-\cT$.
From each subtree $\cT_i$ choose a vertex $t_i$ and collapse each $\cT_i$ in $\cP$ to $t_i$ to obtain a quotient tree $\cP'$ of $\cP$. If an edge in $\cP'$ is labeled with a vertex $t_i'\ne t_i$ from $\cT_i$, then relabel that edge with $t_i$. This turns $\cP'$ into a labeled oriented tree. We say that {\em $\cP$ is injective relative to $\cT$} if $\cP'$ is injective. Note that $\cP$ is injective relative to $\cT$ if and only if: 1) every vertex of $\cP-\cT$ occurs at most once as an edge label in $\cP-\cT$, and 2) every $\cT_i$ contains at most one vertex that is an edge label in $\cP-\cT$. It is clear that if $\cP$ is injective itself, then $\cP$ is also injective relative to $\cT$. We note the following simple and useful observation.

\lemanf\label{W+collapse} If we collapse each connected component $W^+({\bf T}_i)$ of $W^+({\bf T})$ in $W^+(P\cup{\bf T})$ to the vertex $t_i^+$, then we obtain $W^+(P')$. The same is true if we replace $W^+$ by $W^-$.
\lemende
\begin{figure}[ht]\centering
\begin{tikzpicture}
\node at (-2,0) {$\mathcal P$};
\node[red] at (2.5,0.75) {$\mathcal T$};
\fill (0,0) circle (2pt);
\node[below] at (0,0) {$a$};
\fill[red] (1,0) circle (2pt);
\node[red,below] at (1,0) {$b$};
\fill[red] (2,0) circle (2pt);
\node[red,below] at (2,0) {$d$};
\fill[red] (3,0) circle (2pt);
\node[red,below] at (3,0) {$e$};
\fill[red] (4,0) circle (2pt);
\node[red,below] at (4,0) {$f$};
\fill (5,0) circle (2pt);
\node[below] at (5,0) {$c$};

\begin{scope}[decoration={markings, mark=at position 0.5 with {\arrow{>}}}]
\draw [postaction={decorate}] (0,0) -- (1,0) node[midway, above]{$c$};
\draw [red,postaction={decorate}] (1,0) -- (2,0) node[midway, above]{$e$};
\draw [red,postaction={decorate}] (2,0) -- (3,0) node[midway, above]{$f$};
\draw [red,postaction={decorate}] (3,0) -- (4,0) node[midway, above]{$d$};
\end{scope}

\begin{scope}[decoration={markings, mark=at position 0.5 with {\arrow{<}}}]
\draw [postaction={decorate}] (4,0) -- (5,0) node[midway, above]{$a$};
\end{scope}

\node at (-2,-1) {$\mathcal P'$};

\fill (0,-1) circle (2pt);
\node[below] at (0,-1) {$a$};
\fill (1,-1) circle (2pt);
\node[below] at (1,-1) {$b$};
\fill (2,-1) circle (2pt);
\node[below] at (2,-1) {$c$};

\begin{scope}[decoration={markings, mark=at position 0.5 with {\arrow{>}}}]
\draw [postaction={decorate}] (0,-1) -- (1,-1) node[midway, above]{$c$};
\end{scope}

\begin{scope}[decoration={markings, mark=at position 0.5 with {\arrow{<}}}]
\draw [postaction={decorate}] (1,-1) -- (2,-1) node[midway, above]{$a$};
\end{scope}

\begin{scope}[yshift=-3cm]
\node at (-2,0) {$\mathcal Q$};
\node[red] at (2.5,0.75) {$\mathcal T$};
\fill (0,0) circle (2pt);
\node[below] at (0,0) {$a$};
\fill[red] (1,0) circle (2pt);
\node[red,below] at (1,0) {$b$};
\fill[red] (2,0) circle (2pt);
\node[red,below] at (2,0) {$d$};
\fill[red] (3,0) circle (2pt);
\node[red,below] at (3,0) {$e$};
\fill[red] (4,0) circle (2pt);
\node[red,below] at (4,0) {$f$};
\fill (5,0) circle (2pt);
\node[below] at (5,0) {$c$};

\begin{scope}[decoration={markings, mark=at position 0.5 with {\arrow{>}}}]
\draw [postaction={decorate}] (0,0) -- (1,0) node[midway, above]{$c$};
\draw [red,postaction={decorate}] (1,0) -- (2,0) node[midway, above]{$e$};
\draw [red,postaction={decorate}] (2,0) -- (3,0) node[midway, above]{$f$};
\draw [red,postaction={decorate}] (3,0) -- (4,0) node[midway, above]{$d$};
\end{scope}

\begin{scope}[decoration={markings, mark=at position 0.5 with {\arrow{>}}}]
\draw [postaction={decorate}] (4,0) -- (5,0) node[midway, above]{$a$};
\end{scope}

\node at (-2,-1) {$\mathcal Q'$};

\fill (0,-1) circle (2pt);
\node[below] at (0,-1) {$a$};
\fill (1,-1) circle (2pt);
\node[below] at (1,-1) {$b$};
\fill (2,-1) circle (2pt);
\node[below] at (2,-1) {$c$};

\begin{scope}[decoration={markings, mark=at position 0.5 with {\arrow{>}}}]
\draw [postaction={decorate}] (0,-1) -- (1,-1) node[midway, above]{$c$};
\draw [postaction={decorate}] (1,-1) -- (2,-1) node[midway, above]{$a$};
\end{scope}

\end{scope}

\end{tikzpicture}
\caption{\label{reorientexample2} Reorienting in the presence of a sub-LOT.}
\end{figure}
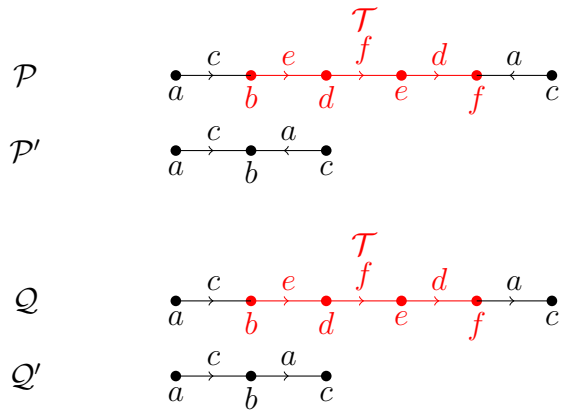

\bspanf\label{secondexample} Figure \ref{reorientexample2} shows labeled oriented trees $\cP$, $\cP'$, $\cQ$, $\cQ'$, and the sub-LOT $\cT$. Let $P$, $P'$, $Q$, $Q'$, and $T$ be the associated LOT-presentations, respectively. Note that $\cP$ is reduced and injective, and the sub-LOT $\cT$ is not boundary reduced (the boundary vertex $b$ does not occur as an edge label in $\cT$). Note further that $\cP'$ is obtained from $\cP$ by collapsing $\cT$ to the vertex $b$, and that $\cQ'$ is obtained from $\cQ$ by collapsing $\cT$ to the vertex $b$. The graph $W^+(P\cup {\bf T})$ is not a tree relative to $W^+({\bf T})$, because if we collapse the connected subgraph $W^+({\bf T})$ of $W^+(P\cup{\bf T})$ to the vertex $b^+$ we obtain the graph $W^+(P')$ (see Lemma \ref{W+collapse}), which is not a tree. It contains a 2-cycle. However we know that we can reorient $\cP'$ to $\cQ'$, so that both $W^-(Q')$ and $W^+(Q')$ are trees. Note that both $\cP'$ and $\cQ'$ already featured in Example \ref{firstexample}. We now reorient $\cP$ outside of $\cT$ to $\cQ$. Both $W^{-}(Q\cup{\bf T})$ and $W^{+}(Q\cup{\bf T})$ are trees relative to $W^{-}({\bf T})$ and $W^{+}({\bf T})$, respectively. If we collapse the connected subgraph $W^{-}({\bf T})$ of $W^{-}(Q\cup{\bf T})$ to the vertex $b^{-}$ we obtain the graphs $W^{-}(Q')$, which is a tree; and if we collapse  the connected subgraph $W^{+}({\bf T})$ of $W^{+}(Q\cup{\bf T})$ to the vertex $b^{+}$ we obtain the graphs $W^{+}(Q')$, which is also a tree.
\bspende

Here is the relative version of Theorem \ref{sprgr}.

\koranf\label{orientationcorollary} If $\cP$ is reduced  and injective relative to $\cT$, then there exists a reorientation $\mathcal Q$ of $\mathcal P$, where only certain edges of $\mathcal P-\mathcal T$ are reoriented, so that $W^+(Q\cup{\bf T})$ and  $W^-(Q\cup{\bf T})$ are trees relative to $W^+({\bf T})$ and  $W^-({\bf T})$, respectively. In fact, $Q$ satisfies the Stallings test relative to $T$. 
\korende

\bewanf
Consider the quotient LOT $\cP'$ obtained from $\cP$ by collapsing each $\cT_i$ to a vertex $t_i$. The LOT $\cP'$ is compressed. In order to see this assume $\cP'$ is not compressed. Then there exists a $\cT_i$ and an edge $e$ in $\cP$ not contained in $\cT_i$ but connected to $\cT_i$ that is labeled with a vertex  $t$ from $\cT_i$. 
Since we assume $\cT_i$ to be a maximal proper sub-LOT the union 
$\cT_i\cup e=\cP$, otherwise $\cT_i\cup e$ would be a larger proper 
sub-LOT. Since $\cT_i$ is a sub-LOT, the vertex of $e$ that is not in 
$\cT_i$ does not occur as edge label in $\cP$. But then $\cP$ is not 
boundary reduced, contradicting our assumption that $\cP$ is reduced.

Theorem \ref{sprgr} implies that there is a reorientation ${\cal Q}'$ of ${\cal P'}$, such that $W^+(Q')$ and $W^-(Q')$ of the LOT-presentation $Q'$ of ${\cal Q}'$ are trees. Let $\cal Q$ be a reorientation of $\cal P$, where 
no edge of $\cT$ is reoriented (so $\cT$ is contained in $\cQ$), so that collapsing each $\cT_i$ in $\cQ$ to the vertex $t_i$ results in $\cQ'$. Let $P$, $P'$, $Q$, $Q'$, $T$, and $T_i$ be LOT-presentations associated with $\cP$, ${\cal P}'$, ${\cal Q}$, ${\cal Q}'$, $\cT$, and ${\cT}_i$, respectively. Collapsing the components $W^+({\bf T}_i)$ of $W^+({\bf T})$ in $W^+(Q\cup{\bf T})$ to $t_i^+$ yields the tree $W^+(Q')$ (see Lemma \ref{W+collapse}). So it follows from Lemma \ref{graphlemma} that $W^+(Q\cup {\bf T})$ is a forest relative to $W^+({\bf T})$. Since $W^+(Q')$ is connected, so is $W^+(Q\cup {\bf T})$. Hence $W^+(Q\cup {\bf T})$ is a tree relative to $W^+({\bf T})$.
In the same way we can argue that $W^-(Q\cup {\bf T})$ is a tree relative to $W^-({\bf T})$.\\

This shows that the forest condition (e) of Definition \ref{relStdef} (the relative Stallings test) holds. The relator conditions (a)-(d) also hold because $Q$ is a LOT-presentation that comes from the reduced LOT $\cQ$. So all relators of $Q$ are cyclically reduced and have exponent sum zero. This implies that all relators of $Q\cup {\bf T}$ have exponent sum zero and hence any word in the generators of $Q$ that represents the trivial element in the group defined by $Q\cup {\bf T}$ has exponent sum zero.\bewende

Let ${\bf x}$ be a set and $w$ be a word in ${\bf x}^{\pm 1}$. Let $S$ be a subset of ${\bf x}$. Define $w_S$ to be the word obtained from $w$ by replacing  $x^{\epsilon}$ in $w$ by $x^{-\epsilon}$, $\epsilon=\pm 1$, 
if and only if $x\in S$. If ${\bf w}$ is a set of words in ${\bf x}^{\pm 1}$, then let ${\bf w}_S$ be the set of words $w_S$, $w\in {\bf w}$. If $P=\langle {\bf x} \ |\ {\bf r} \rangle$ is a presentation, denote by $P_S=\langle {\bf x} \ |\ {\bf r}_S \rangle$.

The map $x\to x^{\epsilon}$, where $x$ is a generator in $P$ and $\epsilon=1$ if $x$ is not in $S$ and $\epsilon=-1$ if $x$ is in $S$, results in a homeomorphism $\phi\colon K(P)\to K(P_S)$ on the corresponding standard 2-complexes. Furthermore, if $f\colon C\to K(P)$ is a vertex reduced spherical diagram, then so is $\phi \circ f\colon C\to K(P_S)$. In particular both $P$ and $P_S$ present the same group. So if $P$ is a LOT-presentation, then $P_S$ is also a presentation of a LOT-group.\\

Let $\cP$ be a reduced labeled oriented tree and let $P$ be the associated LOT-presentation. Let $S$ be a subset of the generators of $P$. Let $\cQ$ be the reorientation of $\cP$ where exactly those edges are reoriented which have their label in $S$. Let $Q$ be the LOT-presentation associated with $\cQ$.\\

The following lemma is essentially Lemma 5.2 of \cite{HR02}. We include
a proof for the convenience of the reader.

\lemanf\label{equalwgraphs} The Whitehead graphs $W(P_S)$ and $W(Q)$ are equal.
\lemende

\bewanf The graphs $W(P_S)$ and $W(Q)$ have the same vertices because the presentations $P_S$ and $Q$ have the same set of generators. In transforming $P$ to $P_S$ an $x^{\epsilon}$, $x\in S$, $\epsilon=\pm 1$, is replaced by $x^{-\epsilon}$ in all the relators  of $P$. Whereas when passing from $P$ to $Q$ an $x^{\epsilon}$, $\epsilon=\pm 1$, is replaced by $x^{-\epsilon}$ only in the relators $r_e$ of $P$, where $e$ is an edge in $\cal P$ with edge label $x$. Thus we obtain $P_S$ from $Q$ by replacing an $x^{\epsilon}$, $x\in S$, $\epsilon=\pm 1$, by $x^{-\epsilon}$ only in the relators $r_{e}$, where $e$ is an edge in $\cal Q$ that contains $x$ as a vertex.

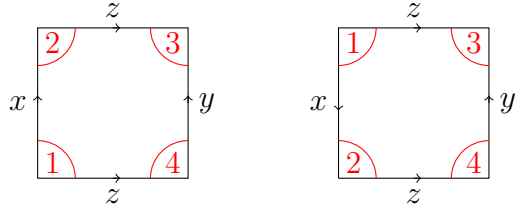
\begin{figure}[ht]
\centering

\begin{tikzpicture}

\draw(0,0)--(2,0)--(2,2)--(0,2)--(0,0);
\draw[->] (1,0)--(1.1,0); \draw[->] (2,1)--(2,1.1); \draw[->] (1,2)--(1.1,2); \draw[->] (0,1)--(0,1.1);
\node [below] at (1,0) {$z$}; \node [above] at (1,2) {$z$}; \node [right] at (2,1) {$y$}; \node [left] at (0,1) {$x$};
\draw [red ] (0.5,0)  arc (0:90:0.5); \draw[red] (0,1.5) arc (-90:0:0.5); \draw[red] (1.5,0) arc (180:90:0.5); \draw[red] (2,1.5) arc (-90:-180:0.5);
\node [red] at (0.2,0.2) {$1$}; \node [red] at (0.2,1.8) {$2$}; \node [red] at (1.8,1.8) {$3$}; \node [red] at (1.8,0.2) {$4$};

\draw (4,0)--(6,0)--(6,2)--(4,2)--(4,0);
\draw[->] (5,0)--(5.1,0); \draw[->] (5,2)--(5.1,2); \draw[->] (4,1)--(4,0.9); \draw[->] (6,1)--(6,1.1); 
\node [below] at (5,0) {$z$}; \node [above] at (5,2) {$z$}; \node [right] at (6,1) {$y$}; \node [left] at (4,1) {$x$};
\draw [red ] (4.5,0)  arc (0:90:0.5); \draw[red] (4,1.5) arc (-90:0:0.5); \draw[red] (5.5,0) arc (180:90:0.5); \draw[red] (6,1.5) arc (-90:-180:0.5);
\node [red] at (4.2,0.2) {$2$}; \node [red] at (4.2,1.8) {$1$}; \node [red] at (5.8,1.8) {$3$}; \node [red] at (5.8,0.2) {$4$};

\end{tikzpicture}
\caption{\label{PtoQ}The relators $r_e=xzy^{-1}z^{-1}$ and $(r_e)_S=x^{-1}zy^{-1}z^{-1}$ contribute the same edges to the Whitehead graphs.}
\end{figure}

Suppose $e=[x,z,y]$ is an edge in $\cal Q$, so $r_{e}=xz(zy)^{-1}=xzy^{-1}z^{-1}$. If neither $x$ nor $y$ are in $S$, then $r_{e}=(r_e)_S$ is also a relator in $P_S$. If $x\in S$ and $y\notin S$, then $(r_e)_S= x^{-1}zy^{-1}z^{-1}$ is a relator in $P_S$. Note that $r_e$ and $(r_e)_S$ contribute the same edges to the Whitehead graph. See Figure \ref{PtoQ}. 
The other two cases, $x\notin S$ but $y\in S$, and both $x$ and $y$ in $S$, lead to relators $(r_e)_S=xzyz^{-1}$ and $(r_e)_S=x^{-1}zyz^{-1}$ in $P_S$ that contribute the same edges to the Whitehead graph as $r_e$. This shows that the Whitehead graphs of $Q$ and $P_S$ are the same.
%%%%%%%%%%%%%%%%%%%%%%%%%%%%%%%%%%%%%%%%%%%%
%%%%%%%%%%%%%%%%%%%%%%%%%%%%%%%%%%%%%%%%%%%%
%%%%%%%%%%%%%%%%%%%%%%%%%%%%%%%%%%%%%%%
%%%%%%%%%%%%%%%%%%%%%%%%%%%%%%%%%%%%%%%
\bewende

Let $\cP$ and $\cT$ be as in the assumption for this section.

\satzanf\label{P_Ssatz} If $\cP$ is injective relative to $\cT$, then there exists a subset $S$ of the generators in $P$ so that $P_S$ satisfies the Stallings test relative to $T_S$.
\satzende

\bewanf By Corollary \ref{orientationcorollary} there exists a reorientation $\cQ$ of $\cP$ so that both
  $W^{\pm}(Q\cup {\bf T})$ are trees relative to $W^{\pm}({\bf T})$. Only edges in $\cP-\cT$ change orientation. Let $S_0$ be the subset of generators of $P$ that occur as labels on edges that change orientation when passing from $\cP$ to $\cQ$. We enlarge $S_0$ to a set $S=S_0\cup \bigcup_{i=1}^nS_i$, where $S_i$ is the empty set if no generator of $T_i$ is contained in $S_0$, and $S_i$ is the set of generators of $T_i$ if $S_0$ contains a generator from $T_i$. Let $\cQ^*$ be the reorientation of $\cP$ where exactly the edges with labels in $S$ are reoriented. Note that $\cQ^*$ contains a reorientation $\cT^*$ of $\cT$, and that $\cQ-\cT=\cQ^*-\cT^*$. Since ${\bf T}={\bf T^*}$ this implies that $Q\cup{\bf T}=Q^*\cup{\bf T^*}$. Thus $W(Q\cup {\bf T})=W(Q^*\cup {\bf T^*})$, and so both $W^{\pm}(Q^*\cup {\bf T^*})$ are trees relative to $W^{\pm}({\bf T^*})$. 

We first show that the forest condition (e) of Definition \ref{relStdef} holds for the pair $T_S\subseteq P_S$. Note that $$P\cup {\bf T}=\langle {\bf x}\  | \ {\bf r}\cup\bigcup_{i=1}^n {\bf u}_i \rangle,$$ where $P=\langle {\bf x}\  | \ {\bf r}\rangle$ and ${\bf u}_i$ is the set of words of exponent sum zero in the generators of $T_i$. Thus 
$$P_S\cup{\bf T_S}=\langle {\bf x}\  | \ {\bf r}_S\cup\bigcup_{i=1}^n {\bf v}_i \rangle,$$ 
where ${\bf v}_i$ is the set of words of exponent sum zero in the generators of $(T_i)_S$. Since the generating sets of $T_i$ and $(T_i)_S$ are the same, we have that ${\bf v}_i={\bf u}_i$ for all $i$. So $$P_S\cup{\bf T_S}=\langle {\bf x}\  | \ {\bf r}_S\cup\bigcup_{i=1}^n {\bf u}_i \rangle.$$ 
Now $$Q^*\cup{\bf T^*}=\langle {\bf x}\  | \ {\bf z}\cup\bigcup_{i=1}^n {\bf w}_i \rangle,$$  where $Q^*=\langle {\bf x}\  | \ {\bf z}\rangle$ and ${\bf w}_i$ is the set of words of exponent sum zero in the generators of $T^*_i$.
Since the generating sets of $T_i$ and $T^*_i$ are the same, we have that ${\bf w}_i={\bf u}_i$ for all $i$. So  $$Q^*\cup{\bf T^*}=\langle {\bf x}\  | \ {\bf z}\cup\bigcup_{i=1}^n {\bf u}_i \rangle.$$ Since $W(P_S)=W(Q^*)$ by Lemma \ref{equalwgraphs} and the same set of relators ${\bf u}=\bigcup_{i=1}^n{\bf u}_i$ is added when enlarging $P_S$ to $P_S\cup{\bf T_S}$ as when enlarging $Q^*$ to $Q^*\cup{\bf T^*}$, it is clear that $W(P_S\cup{\bf T_S})=W(Q^*\cup{\bf T^*})$. Since both $W^{\pm}(Q^*\cup{\bf T^*})$ are trees relative to $W^{\pm}({\bf T^*})$, and $W({\bf T^*})=W({\bf T_S})$, it follows that both $W^{\pm}(P_S\cup{\bf T_S})$ are trees relative to $W^{\pm}({\bf T_S})$. Thus the forest condition holds.\\

We have to check the relator conditions (a)-(d).\\

(a) The relators of $P_S-T_S$ are cyclically reduced because the relators of $P-T$ are cyclically reduced.\\

(b) Note that every relator of $P-T$ is of the form $r_e=ab(bc)^{-1}$, where  $e=[a,b,c]$ is an edge of $\cP$. So every relator contains some generator and its inverse. Hence so does every relator in $P_S-T_S$. Thus no relator of $P_S-T_S$ is a positive or negative word.\\

(c) Let $u$ be a relator of $T_S$. Then $u_S$ is a relator of $T$. Since $T$ is a disjoint union $T=T_1\cup\ldots\cup T_n$, it follows that $u_S$ is a relator in some $T_i$. So $u_S$ is a word in the generators of $T_i$ of exponent sum zero. By construction, the set $S$ contains either no generator of $T_i$, or all of them. Hence $u=(u_S)_S$ is a word in the generators of $(T_i)_S$ of exponent sum zero.\\

(d) Let $w$ be a word in the generators of some $(T_i)_S$ that represents the trivial element in the group defined by $P_S\cup{\bf T_S}=\langle {\bf x}\  | \ {\bf r}_S\cup\bigcup_{i=1}^n {\bf u}_i \rangle$, where $P=\langle {\bf x}\  | \ {\bf r}\rangle$ and ${\bf u}_i$ is the set of words of exponent sum zero in the generators of $T_i$ (see above). Then $w_S$ is a word in the generators of $T_i$ that represents the trivial element in the group defined by $\langle {\bf x}\  | \ {\bf r}\cup\bigcup_{i=1}^n {{\bf u}_i}_S \rangle$. Since by construction $S$ contains either no generator of $T_i$ or all of them, it follows that ${{\bf u}_i}_S={\bf u}_i$ for all $i$.  So $w_S$ is a word in the generators of $T_i$ that represents the trivial element in the group defined by $P\cup {\bf T}=\langle {\bf x}\  | \ {\bf r}\cup\bigcup_{i=1}^n {{\bf u}_i} \rangle$.  Since all relators of $P\cup {\bf T}$ have exponent sum zero, $w_S$ has exponent sum zero. Thus $w_S\in {\bf u}_i$ and hence $w\in {{\bf u}_i}_S={\bf u}_i$. In particular $w$ has exponent sum zero.
\bewende

\section{\label{mainproof}Proof of Theorem \ref{hauptsatz}}

This section is devoted to proving the following theorem.

\satzanf{\label{hauptsatz2}} An
injective labeled oriented tree $\cP$ is aspherical. Furthermore, if $\cH$ is a sub-LOT of $\cP$, then the inclusion induced homomorphism $\pi_1(K(\cH))\to \pi_1(K(\cP))$ is injective.
\satzende

\bewanf We proceed by induction on the number of vertices. If ${\cal P}$ consists of a single vertex the result is true. 

If $\cP$ is not reduced we transform it into a reduced injective labeled oriented tree ${\cP}_{red}$ that contains fewer vertices than $\cP$. Thus, by induction hypothesis, ${\cP}_{red}$ is aspherical, and hence so is $\cP$, because $K(\cP)$ is homotopically equivalent to $K({\cP}_{red})$ (Lemma \ref{LOTtransformation}). Let ${\cal H}_0$ be the image of ${\cal H}$ in ${\cP}_{red}$ under the transformation. By induction hypothesis we know that the inclusion induced homomorphism $\pi_1(K({\cH}_0))\to \pi_1(K({\cP}_{red}))$ is injective. It follows from the commutative diagram in Lemma \ref{LOTtransformation} that $\pi_1(K(\cH))\to \pi_1(K(\cP))$ is injective as well.

From now on we assume that $\cP$ is reduced. Let $\{ {\cal T}_1, \ldots , {\cal T}_n \}$ be the set of maximal proper sub-LOTs of ${\cal P}$. If this set is empty  then $K(\cP)$ is DR and hence aspherical by Theorem 1.1 of \cite{HR01}. Otherwise note that every ${\cal T}_i$ is
compressed and injective and contains fewer vertices than ${\cal P}$. Hence, by induction, each ${\cal T}_i$ is aspherical and $\pi_1(K(\cH))\to \pi_1(K(\cT_i))$ is injective for all sub-LOTs $\cH$ of $\cT_i$.\\

{\bf Case 1.} Suppose that for some $i,j$ we have ${\cal T}_i\cap {\cal T}_j\ne \emptyset$. \\

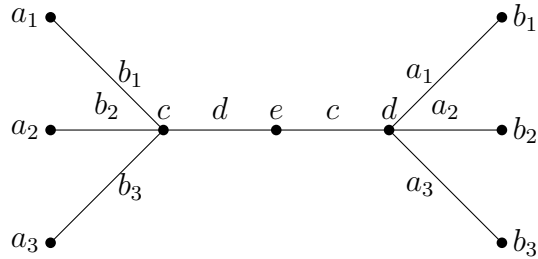
\begin{figure}[ht]\centering
\begin{tikzpicture}
\draw (0,0)--(1.5,0)--(3,0)--(4.5,0)--(6,0); 
\draw (4.5,0)--(6,1.5); \draw (4.5,0)--(6, -1.5);
\draw (1.5,0)--(0,1.5); \draw (1.5,0)--(0,-1.5);
\fill (0,0) circle (2pt); \fill (1.5,0) circle (2pt); \fill (3,0) circle (2pt); \fill (4.5,0) circle (2pt); \fill (6,0) circle (2pt);
\fill (6,1.5) circle (2pt); \fill (6,-1.5) circle (2pt); \fill (0,1.5) circle (2pt); \fill (0,-1.5) circle (2pt);
\node [left] at (0,1.5) {$a_1$}; \node [left] at (0,0) {$a_2$}; \node [left] at (0,-1.5) {$a_3$};
\node [right] at (6,1.5) {$b_1$}; \node [right] at (6,0) {$b_2$}; \node [right] at (6,-1.5) {$b_3$};
\node [above] at (1.5,0) {$c$}; \node [above] at (3,0) {$e$}; \node [above] at (4.5,0) {$d$};
\node [above] at (2.25,0) {$d$}; \node [above] at (3.75,0) {$c$};
\node [right] at (0.75,0.75) {$b_1$}; \node [above] at (0.75,0) {$b_2$}; \node [right] at (0.75,-0.75) {$b_3$};
\node [left] at (4.5+0.75,0.75) {$a_1$}; \node [above] at (4.5+0.75,0) {$a_2$}; \node [left] at (4.5+0.75,-0.75) {$a_3$};
\end{tikzpicture}
\caption{\label{ExCase1} A labeled oriented tree that is the union of three pairwise intersecting maximal sub-LOTs.}
\end{figure}

As an example\footnote{We thank Manuela Ana Cerdeiro for pointing this example out to us.} consider the injective LOT $\cal P$ (with any orientation of its edges) with non-empty intersection of maximal
sub-LOTs as depicted in Figure \ref{ExCase1}. Here let ${\cal T}_i$ be all of $\cal P$ without edges and vertices labeled by 
$a_i$ and $b_i$ for $i=1,2,3$. Then $\{ {\cal T}_1, {\cal T}_2, {\cal T}_3\}$
is the set of maximal proper sub-LOTs of $\cal P$.
%%%%%%%%%%%%%%%%%%%%%%%%%%%%%%%%%%%%%%%
%%%%%%%%%%%%%%%%%%%%%%%%%%%%%%%%%%%%%%%
%%%%%%%%%%%%%%%%%%%%%%%%%%%%%%%%%%%%%%%%
%%%%%%%%%%%%%%%%%%%%%%%%%%%%%%%%%%%%%%%%

We continue with the proof of the theorem in Case 1. We assume
without loss of generality that ${\cal T}_1\cap {\cal T}_2\ne \emptyset$. Then, by maximality it follows that $\cP=\cT_1\cup \cT_2$. Note that there might
be more than two $\cT_i$ (as in the example shown in Figure \ref{ExCase1}).
	The intersection $\cT_{12}={\cal T}_1\cap {\cal T}_2$ is a sub-LOT. Indeed, if $b$ is an edge label in $\cT_{12}$, then $b$ has to be a vertex of ${\cal T}_1$, because ${\cal T}_1$ is a sub-LOT, and $b$ has to be a vertex of ${\cal T}_2$, because ${\cal T}_2$ is a sub-LOT. Thus $b$ is a vertex of $\cT_{12}$. Since $\pi_1(K(\cT_{12}))\to \pi_1(K(\cT_i))$, $i=1,2$, is injective by induction hypothesis, we see that $\pi_1(K(\cP))=\pi_1(K(\cT_1))*_{\pi_1(K(\cT_{12}))}\pi_1(K(\cT_2))$ is an amalgamated product. Furthermore, since both $K(\cT_i)$ and the intersection $K(\cT_{12})$ are aspherical by induction hypothesis, and the inclusion induced maps $\pi_1(K(\cT_{12}))\to \pi_1(K(\cT_i))$, $i=1,2$, are injective,
a theorem of Whitehead \cite{Wh39} (see also Gersten \cite{GerstenBranched}, Theorem 5.1) implies that $K(\cP )$ is aspherical as well.\\

Since $\pi_1(K(\cP))$ is an amalgamated product both inclusion induced homomorphisms $\pi_1(K(\cT_i))\to \pi_1(K(\cP))$, $i=1,2$, are injective. We show that this is true for all proper maximal sub-LOTs and not just for $\cT_1$ and $\cT_2$. Let $\cT_j$ be a maximal proper sub-LOT. If $\cT_1\cap \cT_j=\emptyset$, then $\cT_j\subseteq \cT_2$ and hence $\cT_j=\cT_2$. If $\cT_1\cap \cT_j\ne \emptyset$
and $j\ne 1$, 
then $\cP=\cT_1\cup \cT_j$ and $\pi_1(K(\cP))=\pi_1(K(\cT_1))*_{\pi_1(K(\cT_{1j}))}\pi_1(K(\cT_j))$. In particular $\pi_1(K(\cT_j))\to \pi_1(K(\cP))$ is injective.\\

Now suppose $\cH$ is a sub-LOT of $\cP$. If $\cH=\cP$ then $\pi_1(K(\cH))\to \pi_1(K(\cP))$ is injective. If $\cH$ is a proper sub-LOT, then $\cH$ is contained in some $\cT_j$. Since $\pi_1(K(\cH))\to \pi_1(K(\cT_j))$ is injective by induction hypothesis and $\pi_1(K(\cT_i))\to \pi_1(K(\cP))$ is injective for all $1\le i\le n$, we see that $\pi_1(K(\cH))\to \pi_1(K(\cP))$ is injective.\\

{\bf Case 2.} The ${\cal T}_i$, $i=1,\ldots ,n$, are pairwise disjoint. \\

Let $P$, $T_i$, and $T$ be the LOT-presentations of $\cP$, $\cT_i$, and $\cT$, respectively. By Theorem \ref{P_Ssatz} there exists a subset $S$ of the set of generators of $P$ so that $P_S$ satisfies the Stallings test relative to $T_S$. It follows from Theorem \ref{relativeStallingstest} that $K(P_S)$ is VA relative to $K(T_S)$ and that the inclusion induced maps $\pi_1(K((T_i)_S))\to \pi_1(K(P_S))$ are injective.
Using the homeomorphism $\phi\colon K(P)\to K(P_S)$ defined in Section
\ref{sec:alter} we conclude that $K(P)$ is VA relative to $K(T)$ and the inclusion induced maps $\pi_1(K(T_i))\to \pi_1(K(P))$ are injective. 
By induction hypothesis $\pi_2(K(T_i))=0$ and hence $\pi_2(K(T))=0$. Now $\pi_2(K(P))=0$ follows from Theorem \ref{sVA}. Injectivity of the inclusion induced map $\pi_1(K(\cH))\to \pi_1(K(\cP))$, where $\cH$ is a sub-LOT of $\cP$, follows by the argument given in the paragraph before the statement of Case 2.

This completes the proof of Theorem \ref{hauptsatz2}.
\bewende

\section{Further consequences and applications}\label{sec:consequences}

Many of the arguments given in the previous section can be used to prove asphericity of labeled oriented trees that are not injective. 

\satzanf\label{sgene} Let $\cal P$ be  a reduced
labeled oriented tree and let $\{ {\cal T}_1,\ldots ,{\cal T}_n \}\) be the set of maximal proper sub-LOTs. Let $\cT=\cT_1\cup...\cup \cT_n$. Assume each ${\cal T}_i$ is aspherical.
\begin{enumerate}
\item If the ${\cal T}_i$ are pairwise disjoint and $\cP $ is injective relative to $\cT$, then $\cP$ is aspherical.
\item If ${\cal T}_i\cap {\cal T}_j\ne\emptyset$ for some $i\ne j$, then $\cP=\cT_i\cup \cT_j$.  If both inclusion induced homomorphisms $\pi_1(K({\cal T}_i\cap {\cal T}_j))\to \pi_1(K(\cT_i))$ and 
$\pi_1(K({\cal T}_i\cap {\cal T}_j))\to \pi_1(K(\cT_j))$ are injective, then $\cP$ is aspherical.
\end{enumerate}
\satzende

\bewanf We prove 1:  Let $P$, $T_i$, and $T$ be the LOT-presentations of $\cP$, $\cT_i$, and $\cT$, respectively. By Theorem \ref{P_Ssatz} there exists a subset $S$ of the set of generators of $P$ so that $P_S$ satisfies the Stallings test relative to $T_S$. It follows from Theorem \ref{relativeStallingstest} that $K(P_S)$ is VA relative to $K(T_S)$. Using the homeomorphism $\phi\colon K(P)\to K(P_S)$ defined in Section \ref{sec:alter} we conclude that $K(P)$ is VA relative to $K(T)$. Since we assumed the $\cT_i$ to be aspherical it follows that $\pi_2(K(T))=0$. Now $\pi_2(K(P))=0$ follows from Theorem \ref{sVA}.

The second part of the theorem follows from a theorem of Whitehead \cite{Wh39} concerning the asphericity of unions of aspherical spaces. We used Whiteheads theorem already when we considered Case 1 in the proof of Theorem \ref{hauptsatz2}.
\bewende

The ${\cal T}_i$ may be aspherical for a variety of 
reasons without being injective. For instance they could satisfy small-cancellation
conditions.\\

We conclude this article with an application to long virtual knots. See Kauffman \cite{Kauffman} for an overview of virtual knot theory. A virtual link diagram is a planar 4-regular graph with under- and over crossing information at some nodes.  A virtual knot diagram is a virtual link diagram with only one link component. A long virtual knot diagram  $k$ is obtained by cutting a virtual knot diagram at a point on an edge, thus producing a graph that has exactly two nodes of valency one. A Wirtinger presentation $P(k)$ can be read off in the usual way. It is easy to see that $P(k)=P(\cal P)$, where $\cal P$ is a labeled oriented interval. More details on the connection between labeled oriented intervals and long virtual knots can be found in \cite{HarlanderRosebrock}. We say a long virtual knot diagram is {\em aspherical}, if the standard 2-complex associated with the Wirtinger presentation is aspherical. A virtual knot diagram is {\em alternating} if one encounters over- and under-crossings in an alternating fashion when traveling along the diagram. A long alternating virtual knot diagram is obtained when cutting an alternating virtual knot diagram.

\koranf A long alternating virtual knot diagram $k$ is aspherical. 
\korende

\bewanf The labeled oriented interval that records the Wirtinger presentation of $k$ is injective. The result follows from Theorem \ref{hauptsatz2}.
\bewende

\bigskip\noindent Jens Harlander

\noindent Boise State University, USA

\noindent jensharlander@boisestate.edu

\bigskip\noindent Stephan Rosebrock

\noindent P\"adagogische Hochschule Karlsruhe, Germany

\noindent rosebrock@ph-karlsruhe.de


\begin{thebibliography}{99}

\bibitem{Bogley} W.A. Bogley, {\it J.H.C. Whitehead's asphericity question}, in Two-dimensional Homotopy and Combinatorial Group Theory, eds. C. Hog-Angeloni, W. Metzler and A.J. Sieradski, LMS Lecture Note Series 197, CUP (1993).

\bibitem{BogleyPride} W.A. Bogley, S. J. Pride, {\it Calculating generators of $\pi_2$}, in Two-dimensional Homotopy and Combinatorial Group Theory, eds. C. Hog-Angeloni, W. Metzler and A.J. Sieradski, LMS Lecture Note Series 197, CUP (1993).

\bibitem{Ger87} S. M. Gersten, {\it Reducible diagrams and equations over groups},
in: Essays in Group Theory (S. M. Gersten editor), MSRI Publications 8 (1987), pp.~15--73.

\bibitem{GerstenBranched} S. M. Gersten, {\it Branched coverings of 2-complexes and diagrammatic reducibility}, Transactions of the AMS, Vol. 303, No. 2, (1987), pp.~689--706.

\bibitem{HarlanderRosebrock} J. Harlander, S. Rosebrock, {\it Generalized knot complements and some aspherical ribbon disc complements}, Journal of Knot Theory and Its Ramifications, Volume 12, no. 7, (2003), pp~947--962.

\bibitem{How83} J. Howie, Some remarks on a problem of
J.H.C. Whitehead, {\it Topology} {\bf 22}, (1983), pp.~475--485.

\bibitem{How85} J. Howie, {\it On the Asphericity of Ribbon Disc Complements}, Trans. AMS
289 (1), (1985), pp. 281--302.

\bibitem{HR02} G. Huck and S. Rosebrock,
Weight Tests and Hyperbolic Groups; in {\it Combinatorial and Geometric
Group Theory}; Edinburgh 1993; Cambridge University Press; London
Mathematical Society Lecture Notes Series 204; Editor: A.~Duncan, N.~Gilbert,
J.~Howie; (1995); pp.~174--183.

\bibitem{HR01} G. Huck and S. Rosebrock, Aspherical labeled oriented 
trees and knots, {\it Proc. of the Edinburgh Math. Society} {\bf 44}, (2001), pp.~285--294.

\bibitem{Kauffman} L.H. Kauffman, {\it Virtual knot theory}, Europ. J. Combinatorics, (1999), pp.~663--691.

\bibitem{Rosebrock2} S. Rosebrock, The Whitehead Conjecture - an overview,
{\it Siberian Electronic Mathematical Reports} {\bf 4} (2007), pp.~440--449.

\bibitem{Sta87} J. Stallings, A graph-theoretic lemma and group-embeddings,
in {\it Proceedings of the Alta Lodge 1984}, (ed. Gersten and Stallings),
Annals of Mathematical Studies, Princeton University Press (1987), pp. 145--155.

\bibitem{Wh39} J. H. C. Whitehead, {\it On the asphericity of regions in a 3-sphere}, Fund. Math. 32, (1939), pp~149--166.

\end{thebibliography}
\end{document}